\documentclass{article}
\usepackage{amsfonts, amsbsy, amsmath, amsthm, amssymb, latexsym, enumerate,verbatim,mathrsfs}
\newtheorem{propo}{Proposition}[section]
\newtheorem{defi}[propo]{Definition}

\newtheorem{lemma}[propo]{Lemma}
\newtheorem{corol}[propo]{Corollary}

\newtheorem{theo}[propo]{Theorem}

\newcommand{\bl}{\begin{lemma}\label}
\newcommand{\el}{\end{lemma}}

\newcommand{\ld}{,\ldots ,}

\newcommand{\ra}{ \rightarrow }

\newcommand{\lan}{ \langle }
\newcommand{\ran}{ \rangle }

\newcommand{\diag}{\mathop{\rm diag}\nolimits}

\newcommand{\Id}{\mathop{\rm Id}\nolimits}

\newcommand{\FF}{\mathbb{ F}}

\newcommand{\al}{\alpha}

\newcommand{\lam}{\lambda }

\newcommand{\Om}{\Omega }

\newcommand{\up}{^{-1}}

\newcommand{\si}{\sigma }

\def\d12{{_{12}}}

\def\acf{{algebraically closed field }}

\def\au{{automorphism }}

\def\ccc{{constituent }}

\def\ei{{eigenvalue }}
\def\eis{{eigenvalues }}

\def\f{{following }}

\newcommand{\med}{\medskip}

\def\ii{{if and only if }}

\def\ir{{irreducible }}
\def\irc{{irreducible, }}
\def\irt{{irreducible. }}
\def\irr{{irreducible representation }}

\def\itf{{It follows that }}

\def\mult{{multiplicity }}

\def\po{{polynomial }}

\def\rep{{representation }}
\def\reps{{representations }}
\def\rept{{representation. }}

\def\2k{2^k}

\newcommand{\bp}{\begin{proof} }
\newcommand{\enp}{\end{proof}}

\date{}

\begin{document}
 
\title{Linear groups of alternating type containing non-scalar elements with all but one eigenvalues of multiplicity $1$}
%\author{}

\author{Alexandre Zalesski}

%Address correspondence to
 %E-mail:   alexandre.zalesski@gmail.com (Alexandre Zalesski)}

%\subjclass

\maketitle
 
 \begin{abstract}
Our main goal is to determine   
the irreducible representations of alternating and symmetric groups and their universal central extensions that contain a non-scalar element with all but one eigenvalues of multiplicity 1.  
%This is equivalent to saying that there is only one diminant weight of mult greater than 1.
\footnote{ keywords: alternating groups, irreducible representations, almost cyclic elements}
\end{abstract}

%File AC-ALTI-2025-01-05-final, Directory New projects  Big comp

\section{Introduction}

This paper completes a long-standing project of describing finite primitive \ir  linear groups $G$ containing an element $g$ of prime power order whose all but one \eis have \mult 1. The ground ring is assumed to be a field, usually algebraically closed, and of arbitrary characteristic $\ell$ coprime to $|g|$.  Earlier contributions to this problem are in works   \cite{DZ11,DZ,DPZ,DPZ2,TZ21,TZ2,z23, z24}; more detailed will be given below. An element with this property is called {\it almost cyclic}. If all \eis are of \mult 1 then $g$ is called {\it cyclic}, the term goes back to 
module theory terminology: a module with one generator is called  cyclic. It is known that a module over a cyclic group $\lan g\ran$ is cyclic \ii the minimal \po of $g$ coincides with the characteristic polynomial,  and if $g$ is diagonalizable then this is equaivalent to saying that all \eis of $g$ are of \mult 1. 

The bulk of the project is the case where $G$ is a quasi-simple group, or more precisely those whose derived group  is quasi-simple. (A quasi-simple group $G$ is, by definition, a finite group such that $G=G'$ and $G'/Z(G')$ is a non-abelian simple group.) The case where the simple factor of $G$ is an alternating group is the only one that was not considered yet.
The main goal of this paper  is to prove the \f theorem. 
%Below $F$ is an algebraically closed field of characteristic $\ell$, $A_n$ is the alternating group of $n\geq 5$ letters and $c.A_n$ is the universal  central extension of 
%$A_n$ with center of order $c$. 
The permutational \rep of $A_n$ of degree $n$ is called {\it natural}, and any {\it non-trivial} composition factor of it is called {\it subnatural}.  %\textcolor{red}
(If $n\neq 6$ then all permutational \reps of $A_n$ of degree $n$  are equivalent, whereas those of $A_6$ partition in two equivalence classes; we call
natural the one sending a 3-cycle to a 3-cycle).
% Another one is the twist of it by certain automorphism of $A_6$.)} 
For $g\in G$ the order of $g$
modulo $Z(G)$ is denoted by $o(g)$. 

\begin{theo}\label{th1} Let $G=c.A_n$, $n\geq 5$, be the universal  central extension of $A_n$, and $\phi:G\ra GL_m(F)$ be a non-trivial \irr  of G. Suppose that F is algebraically closed of characteristic $\ell\neq p$ and 
$\phi(g)$ is almost cyclic for some non-central p-element of G. Then $m\leq n+1$, and either $\phi(Z(G))=\Id $ and $\phi$ is subnatural,  
 or one of the cases in Tables $1.1$ and $2.1$ holds,
  in the latter case $m\leq 8$ and $o(g)\leq 9 $.
\end{theo}

In fact we settle a slightly more general case of groups $G$ where   $G'/Z(G')\cong A_n$, see Lemmas \ref{sc0}, \ref{bp4} and \ref{a12}  below.

The assumption that $g\in G$ is of prime power order is not strictly necessary but obtaining a similar result for arbitrary elements of $G$ requires a significally larger work. 
Elements of prime power, in particular,  of prime order are of interest for some applications, see \cite{GPPS,BP,GNPZ, Li}.

  In \cite{GPPS} the authors stated the problem of describing the irreducible subgroups $G$ of  $GL(V)$ (for $V$ being a vector space of finite dimension over a finite field)  that 
contain elements $g\in G$ of prime order $p\not| \, q$ such that  $g$ acts irreducibly on $V/C_V(g)$ (equivalently, on $(\Id -g)V$). They obtained a full description of such groups assuming  that $\dim C_V(g)$ $<\dim V/2$.  Such element $g\in G$ becomes almost cyclic under  a suitable field extension. The results of  \cite{GPPS} were refined in \cite{BP} and  further extended in \cite{DMu}, where for a $p$-element $g$ the assumption  was relaxed to  $|g|>\dim V/3$. In \cite{GNPZ} some more details on the groups listed in \cite{DMu} is given, with  focus on the case where  $\dim C_V(g)=\dim V/2$ and identifying the minimal groups
containing elements in question. Our results together with  \cite{DZ11,DZ, DPZ, DPZ2,TZ2,z24} potentially allows to obtain similar results for subgroups of  $GL(V)$
with no restriction to $\dim C_V(g)$.  Note that the case of   characterisitic 0 ground field is not excluded. Observe that the \ir linear quasi-simple finite groups  over the complex numbers containing a cyclic element have been determined in \cite[Theorem 3.1.5]{KaT1} and a special case of alternating groups has been earlier settled in \cite[Theorem 6.2]{KaT}. In linear group theory  a significant role is played by finite \ir groups generated by reflections and quasi-reflections. A quasi-reflection is a matrix $g\in GL_n(F)$ such that $g-\al.\Id$ is of rank 1 for some $\al\in F$, if $\al=1$ then $g$ is called a reflection. These were classified nearly 50 years ago, see \cite{ZS76,W}. Note that  a quasi-reflection is a special case of an almost cyclic matrix.

 In general, almost cyclic matrices   can be defined as follows \cite[Definition 2.1]{DZ}: 

\begin{defi} Let $M$ be an $(m\times m)$-matrix over a field $F$.
We say that $M$ is almost cyclic  if there exists $\al\in F$ such
that $M$ is similar to $\diag(\al\cdot \Id_k, M_1)$ where $M_1$ is cyclic and $0\leq k<m$.
(Cyclic matrices are exactly those  whose  characteristic \po
coincides with the minimum one.)\end{defi}

In this definition $M$ is not required to be diagonalizable; indeed, if $M$ is unipotent then the Jordan normal form of $g$ is $\diag(\Id_{m-k}, J_k)$, where $J_k$ is a Jordan block of size $k$. Note that the \ir subgroups of $GL_n(q)$ containing a unipotent almost cyclic matrix has been described in \cite{Cra}. 

Finally we observe that if we  assume $\phi(g)$ is cyclic then the conclusion of Theorem 1.1
can be refined by saying that either $\phi$ is subnatural and $o(g)\in\{n,n-1\}$ or
one of the cases of Tables 1,2 holds with $1$ in the columns headed $m$.   

\medskip
{\it Notation} $F$ denotes an \acf of characteristic $\ell\geq 0$ and $\FF_q$ the finite field of $q$ elements. By 
$\mathbb{C}$ and $\mathbb{Q}$ we denote the field of complex and rational numbers, respectively.  By $\mathbb{Z}$ and $\mathbb{N}$ we denote the set of integers and natural numbers.    
For a real number $x\geq  0$  the symbol $\lfloor x\rfloor$ stands for the largest integer  $k$ such that $  k\leq x.$

For a set  $X$ we denote by $|X|$  the cardinality of  $X$ and $Sym(X)$ the group of all permutations of $X$. 

If $G$ is a group then $G'$  is the derived subgroup and  $Z(G)$ the center of $G$.  The identity element of $G$ is usually denoted by $1$. If $g,h\in G$ then $[g,h]:=ghg\up h\up$ and $|g|$ is the order of $g$. If $S\subset G$ is a subset then $C_G(S)=\{g\in G: [x,S]=1\}$  and $N_G(S)=\{g\in G: gS g\up=S\}$. For a subset $X$ of $G$ we denote by $\lan X\ran$ the subgroup generated by $X$. 
$A_n$ and $S_n$ denote the alternating and symmetric groups on $n$ letters.  $c.A_n$, $n>4$  is  a central extension of $A_n$ with center of order $c$ in which every proper normal subgroup lies in its center. 
  $GL_n(F)$ is the group of non-singular $(n\times n)$-matrices over $F$, and we often write $GL_n(F)=GL(V)$ in order to say that $V$ is the underlying vector space for $GL_n(F)$.
%Notations for other classical groups are standard, we usually follow those in \cite{Atlas}, as well as for simple groups.  

If $g\in GL_n(F)=GL(V)$ then $\deg g$ denotes the degree of the minimal \po of $g$.
We say that $g$ is fixed point free if 1 is not an \ei of $g$. The \ei 1 eigenspace of $g$
is often introduced as $C_V(g)$.  The notion of a cyclic and almost cyclic matrix is defined in the introduction. We write $\diag(g_1\ld g_k)$ for a block-diagonal matrix with  diagonal
blocks $g_1\ld g_k$ (not necessarily of the same size). 

Let $\rho:G\ra GL_n(F)$ be a \rept  If $H$ is a subgroup of $G$ then the restriction of $\rho$ to $H$ is denoted by $\rho|_H$. 
%If $G=A_n$, $n>4$, then an   \irr of $G$ of degree at most $n-1$
%(in particular, of degree 1) is called subnatural; more generally, if $G$ is a group such that $G'/Z(G')\cong A_n, n>4$, then an \irr  $\rho$ of $G$ is called subnatural if $\rho(Z(G'))=\Id$ and   $\dim \rho<n$. 

\section{Preliminaries}

\def\sft{Suppose first that }

If an almost cyclic matrix $M$ acts in a vector space $W$ and $U$
is an $M$-stable subspace of $W$ then the restrictions of $M$ to
$U$ and to $W/U$ are almost cyclic matrices.

\bl{va1}  Let $g\in G\subset GL_n(F), n>1$ be a finite \ir subgroup 
  generated by $k$ conjugates of $g$. Suppose that g is almost cyclic,  
and $d=\dim(g-\lam\cdot  \Id)V$ for some $\lam\in F$. Then  $n\leq dk$ and $n\leq (o(g)-1)k$. 
\el

\bp Let $g_1\ld g_k$ be conjugates of $g$. Then $G=\lan g_1\ld g_k\ran $ stabilizes the subspace W:=$(g_1-\lam\cdot  \Id)V+(g_2-\lam\cdot  \Id)V+\cdots +(g_k-\lam\cdot  \Id)V$. As $G$ is irreducible, we have $W=V$.  So the first inequality follows. 
For the second see \cite[Lemma 2.11]{DPZ2}. \enp

\begin{lemma}\label{ma2} {\rm \cite[Lemma 3.2]{z23}} Let $M=M_1\otimes M_2$ be a Kronecker product of diagonal non-scalar matrices $M_1, M_2$ of sizes $m\leq n$,
respectively.

$(1)$ The \ei multiplicities of $M$ do not exceed $mt$, where t is the maximal \ei \mult 
of $M_2$. 

$(2)$ If $M_1,M_2$ are cyclic then the \ei multiplicities of $M$ do not exceed  $m$.

$(3)$ Suppose that M is almost cyclic. Then $M_1$ and $M_2$ are cyclic.

$(4)$ Suppose that M is almost cyclic  and $M_i$ is similar to $M_i\up$ for $i=1,2$. Then the \ei multiplicities of $M$ do not exceed  $2$. In addition, if $e$
is an \ei of M of \mult $2$ then $e\in\{\pm 1\}$. \el

\section{Representations of $A_n$ of small degrees }

%In this section we determine  the \ir $F$-representations of the universal covering of the alternating $A_n$ with $n>4$
%containing an almost cyclic matrix.

Let $G$ denote $A_n$ or $S_n$,  $n>4$. %We keep the notion of a natural permutation \rep of $A_n$ for analogues \rep of $S_n$. 
We say that an \irr $\phi$ of $S_n$ is {\it subnatural} if the \ir constituents of $\phi|_{A_n}$ are subnatural \reps of $A_n$. This extends to $S_n$ the notion of a   subnatural \rep of $A_n$ given in the introduction. In fact, any subnatural \rep of $S_n$ is \ir on $A_n$.

In this section we denote by $P_n$  the permutation $FS_n$-module of dimension $n$ arising from 
a natural \rep of $S_n$. It is well known that $P_n$ has a unique non-trivial composition factor, and one or two trivial factors depending of whether $\ell\nmid n$ or $\ell|n$, respectively. 
We denote the non-trivial factor by $W_n$, and denote by $W_n'$  
the restriction of $W_n$ to $A_n$. Then $W_n'$ remains \irt  
%(If $n=6$ then there are 2 non-equivalent permutation representations of $S_6$ of degree 6, we could extend the this notation to declare both of them natural.)
If $\ell\neq 2$  then we set $P_n^-=P_n\otimes \al_n$, where $\al_n$ is a non-trivial one-dimensional $F$-\rep of $S_n$, and  $W_n^-=W_n\otimes \al$.  We call $W_n,W_n^-$ and $W_n'$   {\it standard} $FG$-modules (some authors call them deleted). The representations of $G$ afforded by standard $FG$-modules are   subnatural.
 It is well known that $\dim W_n=\dim W_n'=n-1$ if $\ell \nmid n$, and $n-2$ otherwise.   If $n>8$ then the converse holds, in the sense that every $FG$-module of dimension $d$ with $1<d<n$ is standard. See \cite{W2,W3}.

\def\pp{permutation }

\bl{kz1} {\rm \cite[Corollary 2.4]{kz}} Let $n>7$ and let $\phi$ be a non-trivial  \irr of $A_n$ such that all composition factors of $\phi|_{A_{n-1}}$ are     either trivial or  subnatural \reps of $A_{n-1}$. Then $\phi$ is a subnatural \reps of $A_{n}$.\el

Using induction on $n$ with $n=8$ the induction base, we have: 

\begin{corol}\label{kz2} Let $6<k<n $ and let $\phi$ be a non-trivial  \irr of $A_n$ such that all composition factors of $\phi|_{A_{k}}$ are    either trivial or 
subnatural \reps of $A_{k}$. Then $\phi$ is a subnatural \reps of $A_{n}$.
\end{corol}

 Note that all subgroups of $A_n$,  $n>6$,  isomorphic to $A_{n-1}$ are conjugate.  
(By \cite[\S 5.2]{DM} $A_{n-1}$ is intransitive, and, by an order reason,  a subgroup $X\cong A_{n-1}\subset A_n$ fixes a point  and even coincides with the stabilizer of a point; as $A_n$ is transitive, the point stabilizers are conjugate in $A_n$.) 
 
\bl{vw1}  $(1)$ The alternating group $A_n$ with n odd  is generated by two cycles of order $n$.  
 
$(2)$ For  $n>6$ even the group $A_n$ is generated by two elements of order $n/2$ of double cycle structure, and $S_n$  is generated by two cycles of order $n$. 
\el

\bp (1) are well known. For a proof, observe that every $g\in A_n$ is a product of two $n$-cycles, see for instance  \cite{dw1}. For $n=3$ the statement is trivial, for $n=5$ this follows by inspection of  the subgroups of $A_5$ in \cite{Atlas}. 
Let $n>5$. By the Bertrand postulate, if $n-2>3$ then there exists a prime  $p>2$ such that $(n-2)/2<p< n-2$, see \cite{bp1}.   Let   $g\in A_n$ be  a $p$-cycle %with $|g|=p$ 
and let  $g=xy$ for some $n$-cycles $x,y\in S_n$.  
So $X=\lan x,y\ran$ is a transitive group containing a $p$-cycle with $(n-2)/2<p< n-2$. Then 
$X$ is primitive (otherwise the imprimitivity blocks are of size at most $n/3$ as $n$ is odd, which leads to a contradiction). By  \cite[Theorem 3.3E]{DM},   $A_n\subseteq X$.  

(2). We can assume that $A_n$ acts on $\Om=\{1\ld n\}$ and  $h=(1\ld n/2)((n/2)+1\ld n)$. Let $t=(1,
2, n)$. Then $x:=tht\up h\up=(1,2,(n/2)+1, 3,n)$. Then 
$X=\lan h, tht\up \ran$ is primitive.  (Indeed, suppose the contrary. Then $x\in X$ cannot permute the imprimitivity blocks as $x$ fixes $n-5$ points of $\Om$. So $x$ stabilizes them, and hence $\{1,2,(n/2)+1, 3,n\}$ lie in  
one of them. Denote it by $\Om'$. As $h(1)=2\in\Om'$, we observe that $h(\Om')=\Om'$, and hence the $h$-orbits of $1$ and $n$ are contained in $\Om'$. So $\Om'=\Om$.) 
By  \cite[Theorem 3.3E]{DM}, either  $n\leq 7$ (which is not the case by assumption) or  $A_n=  X$.  

The claim on $S_n$ follows from that for $A_n$; indeed,   if $g\in S_n$ is an $n$-cycle then $g^2\in A_n$ has the double cycle structure.\enp
 
\bl{ku6} Let F be an \acf of characteristic $\ell\geq 0$,  $G$  a finite group and %where $m=p^a$ for a prime p. 
let $\phi:G\ra GL_m(F)$ be an \ir F-\rep of G. Let $g\in G$ and let h  be the projection of $g$ into $G/Z(G)$. Suppose that $\phi(g)$ is  almost cyclic. 

$(1) $ Suppose   $G/Z(G)=A_n$ with n odd, and $h\in A_n$ is an n-cycle.%of order m 
%and     $\phi(g)$ is  almost cyclic. 
Then $m\leq 2 (n-1)$. 

$(2) $ Suppose that $G/Z(G)=S_n$ with n even and  $h\in S_n$ is an n-cycle.   Then $m\leq 2 (n-1)$. 

$(3) $  Suppose that $G/Z(G)=A_n$ with n even and  $h\in A_n$ is a double $(n/2)$-cycle.   Then $m\leq n-2$.  
\el

\bp  By Lemma \ref{vw1}, $G/Z(G)$ is generated by two conjugate of $h$. Let $X$ be a subgroup generated by two conjugate of $g$ whose projection in  $G/Z(G)$ generate $G/Z(G)$. Then $G=Z(G)\circ X$, a central product. Then  (1) and (2) follow from Lemma \ref{va1}. (3) As $|g|=n/2$, we similarly
have  $m\leq 2(|g|-1)=n-2$.\enp

\begin{theo}\label{wo2} {\rm \cite[Theorem 1.3]{W1}} Let $ n = 2^{ w_1}+ 2^{ w_2}+... + 2^{ w_s}$  with $w _1 > w_2 > ... > w_s$. Also let F be any field 
 with characteristic $\neq 2$. Then:

$(1) $ the degree of any   faithful representation of $2. S_n$ over F is divisible by $2^{\lfloor(n-s)/2\rfloor};$

$(2)$ the degree of any  faithful representation of $2. A_n$ over F is divisible by
$2^{\lfloor(n-s-1)/2\rfloor}.$\end{theo}
%$\lfloor r \rfloor$

\bl{ee2} In notation of Theorem {\rm \ref{wo2}} we have $2^{\lfloor(n-s-1)/2\rfloor}>2(n-1) $ for $n>13$.  \el

\bp Suppose the contrary, that $2^{\lfloor(n-s-1)/2\rfloor}\leq 2(n-1) $. Then $2^{(n-s-3)/2}\leq 2(n-1)$ and
$2^{n-s-3}\leq 4(n-1)^2 $,  $2^{n-3}\leq (n-1)^22^{s}$. Note that $1+2+\cdots +2^{s-1}=2^{s}-1\leq n$, so $2^{n-3}\leq (n-1)^2 (n+1)$ and $n\leq 14$.
%or $2^{n-5}\leq (n-1)^2$.
One checks that $2^{\lfloor(n-s-1)/2\rfloor}>2(n-1) $ for $n=14$. \enp

 In Lemma \ref{33e} we write $[k, l]$ for a permutation with two cycles of size $k,l$. 

\bl{33e} Let $G=A_n$ or $S_n$, $n>4$,  and, for a prime p, let $1\neq g\in G$ be a p-element, $|g|=p^a$. Suppose that $p\neq \ell$ and g is almost cyclic on $W_n$.  Then  either  $g=[1^{n-p^a}, p^a]$, or $p=2$,   $n=2^a+2$,   $g=[2,n-2]$.

In addition, $\deg g\geq |g|-1$, and the equality holds \ii $|g|=n$ and g is cyclic, or $p=2$, $n=2^a+2$,   $g=[2,n-2]$ and the \ei $-1$ \mult equals $2$.
 \el

\bp  Observe that $P_n$ is a direct sum of $P_c$, where $c$ runs over the  sizes
of the cycles in the cycle decomposition of $g$. The characterrstic polynomial of a $c$-cycle element in $P_c$ is $x^c-1$, so the \eis of this element
in $P_c$ are all $c$-roots of unity, each occurs  with \mult 1. Let $[c_1^{m_1}\ld c_k^{m_k}]$ be a cycle decomposition of $g$, where we can assume that $c_1<\cdots <c_m$. Then the characteristic polynomial of $g$ in $P_n$ is $\Pi_{i=1}^k (x-c_i)^{m_i}$.  As $g$ is a $p$-element, %we have $c_i|c_j$  and 
$(x-1)|(x^{c_i}-1)$  for $i=1\ld k$.
\itf every non-trivial $p$-root of unity is an \ei of $g$ in $P_n$ of \mult $\sum_{c_i\neq 1} c_i^{m_i}$, whereas the  \mult of \ei 1 equal  $\sum c^{m_i}_i$. Therefore, 
 $g$ is almost cyclic in $P_n$ \ii 
$g=[c_1] $  or $[1^m,c_2]$ for some  $c_1,c_2>1$. 

 Recall that the composition factors of $P_n$ are $W_n$ and $1_G$, the latter appears with \mult 1 if $n$ is not a multiple of $\ell$, and 
with \mult 2 otherwise. Clearly, the multiplicities of all \eis $\lam\neq 1$ of $g\in G$ in $P_n$ and $W_n$ are the same. If $\lam\neq \pm 1$ then this occurs with \mult 1 in 
$P_n$. Therefore,  $g$ has at most one cycle of size greater than 2. If $g$ has no cycle of size 2 then $g=[1^{m_1},c_2]$; this case is recorded in the lemma conclusion.

 Suppose that $g$ has a cycle of size 2 and $g\neq [1^{m_1},2]$.   Then $p=2$ and $-1$  \eis of $g$ on $W_n$ of \mult at least 2. Therefore,  $g=[2,c_2]$ (otherwise each 1 and $-1$ are \eis of $g$ on $W_n$ of \mult at least 2).  Whence the result. \enp

\bl{a12} For $4<n<13$ let $G/Z(G)=A_{n}$. %  or $S_n$. %and $|Z(G)|\leq 2$. 
Let $\phi:G\ra GL_m(F)$ be a faithful  \ir  \rep of G  and $g\in G$ a p-element for $p\neq \ell$.  Suppose that $\phi(g)$ is almost cyclic. Then either $Z(G')=1$ and $\phi|_{G'} $ is subnatural  
 or  $n\leq 10$ and one of the cases in Tables $ 1.1$ and $2.1$  holds. 

%In addition,the cases where $G/Z(G)=S_{n}$ and $\phi|{G'}$ is reducible are listed in Table $3.$
\el
 
\bp In \cite{AB} one finds the Brauer characters of the universal covering of the simple groups $A_{n}$ for $4<n<13$ and their extension by an outer automorphisms. So the result follows by inspection of the Brauer characters of these groups. In particular,      the tables contain no entry   for $n=11,12$.
%the subnatural \irr of $2.A_n$ is the only 
%one where some $p$-element with $2<p\leq n$ is almost cyclic. 
The data of Tables 1.1 and 2.1 are collected from \cite{AB}. \enp

It seems to be useful, for readers' convenience, to  comment  special cases, where our conclusion, for fixed $|g|$ and $m=\dim\phi$, depends on the  conjugacy class of $g$ and  the choice of $\phi$.

Suppose that   $\phi$ is not subnatural.  
  If $n\neq 6,7$ then $G'$ is isomorphic to $A_n$ or $2A_n$, so $g\in G'$ if $|g|$ is odd. 

Let $n=10$. If $\ell=5$ then there are 2 faithful \ir \reps  of $2.A_{10}$ of degree 8. 
If $|g|=9$ then $g$ is cyclic in one of these \reps and almost cyclic  in the other one with \ei 1 of \mult 2.   

 The case with $G=2.A_{10}$ and $o(g)=8$ has a feature that does not occur in other cases in the tables. In this case $G$ has 2 non-equivalent faithful \ir \reps $\phi_1,\phi_2$ of degree 8 and 2 conjugacy classes
of order 8 which glue in $G/Z(G)$. In fact, these are $g$ and $zg$ where $1\neq z\in Z(G)$. In each case $\deg\phi_i(g)=\deg\phi_i(zg)=7$ for $i=1,2$ and the \ei $e$ of \mult 2 is 1 or $-1$. Clearly, for a fixed $i$ we have $e(g)=-e(zg)$. The feature is
$\phi_1(g)=-\phi_2(g)$ and $\phi_1(zg)=-\phi_2(zg)$. This fact is reflected in Table 2.1
by writing $e=\pm 1$ in one \rep and $e=\mp 1$ in the other one.

Let $n=9 $, $o(g)=9$ and $\ell\neq 2,3$. Then $2A_9$ has  two faithful  \ir representations $\phi$ of degree 8 and two conjugacy classes of elements of order 9. % If $\phi$ is subnatural then $\phi(g)$ is cyclic, $\deg\phi(g)=8$ and 1 is not an \ei of $\phi(g)$. Otherwise, 
%Depending on the class and the represenation, 
In one of these two \reps $\phi(g)$ is cyclic  and  1 is not an \ei of $\phi(g)$,  in the other \rep  $\phi(g)$ is almost cyclic (not
cyclic), $\deg\phi(g)=7$ and  1 is  an \ei of $\phi(g)$ of \mult 2.  If $\ell=2$
then there is two non-subnatural \ir \reps of degree 8. As for $\ell>3$, in
one of these two \reps $\phi(g)$ is cyclic  and  1 is not an \ei of $\phi(g)$,  in the other \rep  $\phi(g)$ is almost cyclic (not
cyclic), $\deg\phi(g)=7$ and  1 is  an \ei of $\phi(g)$ of \mult 2. 

%The same is true for $\ell=2$, where there are two   non-subnatural \irr of degree 8. 

% As $\phi(Z(G'))\in \{\pm \Id\}$, %=Z(\phi(G))\circ \phi(2.S_9)$,
% we observe that $\phi(zg)=-\phi(g)$ is almost cyclic, and  \ei $-\sqrt{-1}$ occurs with \mult 2. This is also true for elements $g\in 2.S_8$ with $o(g)=8$.  

%for every element $g\in G$ with $o(g)=8$, and the \ei of \mult 2 can be any $8$-root of  $\lam$. 

Let $n=8$ and $\ell= 2$.   Recall that $A_8\cong SL_4(2) $ and $|g|\in\{3,5,7\}$.  %Suppose  that , and %$|g|$ is odd. 
%Then $|\phi(Z(G))|$ is odd and $\phi(G)=\phi(Z(G))\times A_8$ or $\phi(Z(G))\times S_8$.  In this case $|g|$ is odd, we can assume that $g\in G'\cong A_8$,and hence
 There are two non-equivalent \ir \reps of $G$ of degree 4. One  observes that $\phi(g) $ is cyclic for $|g|=5,7$. There are two classes of 
elements $g$ of order 3, and $\phi(g) $ is almost cyclic with $\deg\phi(g)=3$ exactly for one of two classes.

 Let $n=7$ and $G=6.A_7$. % If $\dim\phi=3$ then $\ell=5$ 
%then %$G$ has an \irr $\phi$ of degree 3, and $\phi(g)$ is cyclic for every $g\in G\setminus Z(G)$ for $o(g)>2$; if $o(g)=2$ then $\phi(g)^2=1$ and $g$ is almost cyclic with $-1$ to be an \ei of \mult 2. 
Then   $G$ has two  \ir \reps $\phi$ of degree 4 for $\ell\neq 7$, $|\phi(Z(G))|\leq 2$, and  $|\phi(Z(G))|=1$ \ii $\ell=2$.  These extend to $2.S_7$ \ii $\ell=7$.
The group $A_7$ has two conjugacy classes of elements of order 3. If $g$ corresponds to a 3-cycle element in $A_7$ then $\phi(g)$ is  almost cyclic  whenever $\dim\phi=4$.
If $g$  corresponds to a double 3-cycle element in $A_7$ then $\phi(g)$ is  not  almost cyclic.

Let $n=6$. If $\dim\phi=4$ and $o(g)=3$ then, as in the case with $n=7$, there are two conjugacy classes of elements of order 3 in $A_6$, and $\phi(g)$  is  almost cyclic in one of them and not almost cyclic in the other (for $\ell\neq 3$). Note that one of these \reps is subnatural, and if $\phi$ is not subnatural then  $\phi(g)$  is  almost cyclic when $g$ corresponds to a double 3-cycle element in $A_6$.  
 
Let $\dim\phi=5$. Then $Z(G)=1$, $\ell\neq 2$ and $G$ has two \irr of degree 5, one of them is subnatural. (These are obtained from each other by a twist with an outer \au of $A_6). $ 
If $o(g)=3$ then there are two conjugacy classes of elements of order 3, elements of  one of them are almost cyclic and $\deg\phi(g)=3$, those in the other class are not  almost cyclic. 
 \med

Note that ${\rm Out}\,(A_n)$ has order 2 if $6\neq n>4 $ and of order 4 (and of exponent 2) for $n=6$.
 If $G'/Z(G')\cong A_n$ and $g\in G$ is of odd prime power order then $g\in G'Z(G)$.  
For $4<n<13$ Lemma \ref{a12} provides a sufficient information on 
almost cyclic elements in irreducible \reps $\phi$ of such groups provided $\phi|{G'}$ is irreducible. Below we complement Lemma \ref{a12} by adding the cases where  $\phi(g)$  is almost cyclic and  $\phi|_{G'}$ is reducible.   

\bl{sn1}  Let $G/Z(G)\cong S_n$, $4<n\leq 13$, and let $\phi$ be an \ir $\ell$-\rep of $G$ such that $\phi|_{G'}$ is reducible. Let $p>2$ be a prime, $p\neq \ell$ and let $g\in G$ be a p-element. Suppose that $\phi(g)$ is almost cyclic. Then one of the cases of Table $1.2$ holds. \el

 The entries of Table $1.2$ are extracted from \cite{AB}.
Note that ${\rm Aut}\,A_6$ have no \ir projective \rep  that is reducible on 
every normal proper subgroup.  
(Otherwise, by Clifford's theorem, the group
$c.A_6$, for some fixed $c$ and $\ell$, would have 4 distinct \ir Brauer characters of the same degree that are conjugate in ${\rm Aut}\,A_6$. This is not the case by \cite{AB}.)

\med
 
Let $p=2$ and $\ell\neq 2$. Table 2.1 is obtained by inspection of \cite{AB}.  In this table $G/Z(G)=A_n$ and we differ the groups in question by indicating $G'$ in the 1st column. 

Table 2.2 lists the cases where $G=c.A_n.2$, a non-split semidirect product of $G=c.A_n$ and a group of order 2. As we shall see, this is sufficient for using in the proof of Theorem 1.1 and describe details on maximum \ei \mult for $g\in G$.

% (an element of   is a proper subgroup of $ G/Z(G)$ and $G=\lan g,G'\ran$. 
Recall that $G/Z(G)\cong S_n$ if $n\neq 6,n>4$, whereas 
%If $n=6$ requires more attention as 
the group  $ A_6$ has 3 outer automorphisms of order 2, and hence 3 non-isomorphic groups $A_n.2$.  Following \cite{AB} we denote the projection of $g$ into $G/Z(G)$ by $2_i$, $i=1,2,3$, where $i=2$ if $G/Z(G)\cong S_n$,  $i=1$ %.  In addition, $A_6\cong PSL_2(9)$ and $2_1$ is the involution in 
if $G/Z(G)\cong PGL_2(9)$, and $2_3=2_1\cdot 2_2$.
 
For every \irr $\phi$ of $G$ one can consider $\phi\otimes\al$, where $\al$
is a non-trivial one-dimensional \rep of $G$. As $G=\lan g,G'\ran$, we have $\al(g)=-1$.
 If $e$ is a unique \ei of $\phi(g)$ of \mult $m>1$ then  $-e$ is an \ei of $(\phi\otimes\al)(g)$ of \mult $m>1$ too. In particular,  $\phi$ and $\phi\otimes\al$ are not equivalent. In order to avoid repetitions of these very similar cases we write $\pm e$ in the column headed $e$. 
  
 The only case in Table 2.2 where $\phi|_{G'}$ is reducible is the one  of degree 4
written at the 3rd row.   
 
\section{The main result}

\bl{ct7} Let $G$ be a group such that $G/Z(G)\cong  A_n$ or $S_n$, $n>9$,  and let $g\in G$ be a p-element for a prime p. Let $\ell\neq p$ and let $\phi: G\ra GL(m,F)$ be a faithful \irr such that $\phi (g)$  is almost cyclic. Suppose that   $n=o(g)$.  Then $\phi|_{G'}$ has at most two \ir constituents, both of them   are subnatural.\el  
 
 \bp  We identify $G/Z(G)$ with a subgroup of $S_n$. 
Let   $h$ be the projection of $g$ into $G/Z(G)$. Then $h$ is a cycle of length $n$ (as $n=o(g)$).    
   By Lemma \ref{vw1}, if $p=2$ then $G/Z(G)\cong S_n$ is generated by two conjugates of $h$; if $p>2$ then $A_n$ is generated by two conjugates of $h$. Let $h'$
be a conjugate of $h$ such that $\lan h,h'\ran\cong S_n$, resp., $A_n$ if $p=2$, respectively, $p>2$. Let $g'\in G$ be an element such that  $h'=gZ(G)$.  
 Then   $\lan g,g'\ran$ projects onto $\lan h,h'\ran$. Set $G_1=\lan Z(G), g,g'\ran$.
Then $G_1=G$ if $p=2$ or $p>2$ and $G/Z(G)=A_n$, otherwise $|G:G_1|=2$ and $G_1/Z(G_1)\cong A_n$.
In the latter case $\phi|_{G_1}$ has at most 2 \ir constituents. Let $\phi_1$ be one of them; then $\dim\phi=a\dim\phi_1$ with $a\leq 2$.  
 
    As $\phi$ is irreducible, the dimension of $ \phi$, resp., of $\phi_1$ does not exceed $2o(g)-2$ by Lemma \ref{va1}.  In each case, the dimension of an \ir constituent $\tau$ of $\phi|_{G'}$ do not exceed $2n-2$.

 Suppose first that $\phi(Z(G'))\neq 1$. Then $\ell\neq 2 $ and $\tau(G')\cong 2.A_n $.  By Lemma \ref {wo2}, 
$\dim \tau \geq 2^{\lfloor(n-s-1)/2\rfloor}$, where $s$ is the number of non-zero terms in the
2-adic expansion of $n$.  So $\dim\tau\leq 2n-2$, which  implies $n\leq 13 $ by Lemma \ref{ee2}.  Let $n=13$. Let $\rho$ be a faithful \irr of $G'=2.A_{13}$. We show that
 $d\geq 32$. Indeed, let $X\cong 2.A_{12}\subset G'$. Then the \ir constituents 
of    $\rho|_X$ are faithful \reps of $X$. By  \cite{AB}, the degree $d'$, say, of a  faithful \ir \rep $\rho'$ of $X$     is not less that 32, unless $\ell=3$
where $d'\geq 16$. Moreover, if $d'>16$ then $d'\geq 144$ for $\ell=3$. So we are left with $\ell=3$, $d'=16$ and $\rho|_X$ is \irt Let $x\in X$ be of order 11. By \cite{AB}, 
the Brauer character of $\rho'$, and hence of $\rho$, takes an irrational value at $x$.   
This is a contradiction as $x$ is conjugate to every $x^i$ for $0<i<10$, so $\rho(x)$ must be an integer. (Note that $|N_{G'}(\lan x \ran)|=10$ as $G'/Z(G')$ contains
$S_{11}$.)

Suppose  that $\phi(Z(G'))=1$ so $\phi(G')=A_n$. By  James \cite[p. 420 and Theorem 7]{J}, if $n>14$ and $\si$ is an \ir modular  \rep of $S_n$ then
either all \ir constituents of $\si|_{G'}$ are subnatural or $\dim\si>n(n-5)/2$. 
If $n>14$ and $\tau$ is not subnatural then $\dim\tau\leq 2n-2$, so $\dim\phi\leq 4n-4$, whence  $n\leq 12$.
(Indeed, we have  $n(n-5)/2\leq 4n-4$ implies $n(n-5)\leq 8n-8$, $n^2-13n+8\leq 0$.)  As $n$ is a prime power, $n\leq 11$ or $n=13$. If $\phi|_{G_1}$ is \ir then $\dim\phi\leq 2n-2$, whence $n(n-5)/2\leq 2n-2$, a contradiction for $n=13$.
If $\phi|_{G_1}$ is reducible then  $\dim\tau=\dim\phi_1\leq 2n-2$, which is 24 for $n=13$. If $\tau$ is not subnatural then $\dim\tau>12$, and then $\dim\tau\geq 32$
by \cite{at13}, a contradiction. So $n\leq 11$.

 As $n$ is a $p$-power,  $n\neq 10$; for  $n= 11$ the lemma is true by 
Tables $1.1$ and $1.2$.
\enp

\bl{sc0} Let $p=2$ and let G be a group such that $G'\cong c.A_n$, $n>11$, and $G=\lan g,G'\ran$, where   $g\in G$ is a $2$-element.  Let $\phi: G\ra GL(m,F)$ be a   faithful \irr of G such that  $\phi (g)$  is almost cyclic. Then  $G'\cong A_n$ 
% $\phi(Z(G'))=\Id$ 
and   $\phi|_{G'}$ is a subnatural.\el 

\bp  Note that $G/Z(G)\cong A_n$ or $S_n$ and $c=|Z(G')|\leq 2$. Let $h$ be the projection of $g$ into $G/Z(G)$, and let $k$ be the largest cycle size in the cycle decomposition of $h$. Note that $k=|h|=o(g)$. 

Suppose first that $o(g)\geq 8$.
 If $k=n$ then the result is contained in Lemma \ref{ct7}. So  $8\leq k<n$. Let $\Om$ be the natural permutational set for $G/Z(G)$, and let  $\Om= \Om_k\cup \Om_{n-k}$, where  $\Om_k$ is an $h$-orbit of size $k$ and $\Om_{n-k}$ is the complement of  $\Om_k$ in  $\Om$. Let $X_1\cong A_k$, $X_2\cong A_{n-k}$  be the groups of all odd permutations on  $ \Om_k,\Om_{n-k}$, respectively, and let $ X=X_1\times X_2\subset Sym(\Om)$.   

If $h\notin X$ then $h^2\in X$ (as  $h$ preserves $\Om_k$ and $ \Om_{n-k}$). 
Let $G_1,G_2$ be the preimages of $X_1,X_{2}$ in $G$. Then   $G_i$   are $g$-invariant for $i=1,2$ and $g^2\in G_1G_2$. In addition, $[G_1,G_2]=1$.  (Indeed,   $X_1$   is perfect, and $[G_1,G_2]\subset Z(G)$. Then  for $y\in G_2$ the mapping $G_1\ra Z(G)$ defined by $x\ra [x,y]$ $(x\in G_1)$ is a homomorphism, which is trivial as $X_1$ is perfect. So $[x,y]=1$.)

 Set $G_3=\lan g,G_1\ran$. Then $G_3'=c.A_k$, where $c\leq 2$,  and $\tau(g)$ is almost cyclic for every \ir constituent $\tau$ of $\phi|_{G_3}$. As $k\geq 8$, the action 
of $g$ on $G_1$ is induced by an inner automorphism of $S_k$. \itf $G_3/Z(G_3)$ is isomorphic to $A_k$ or $S_k$. Since $o(g)=k$ in $G_3$, By Lemma \ref{ct7}, we conclude that $\phi|_{G_3'}$ is either subnatural or trivial.  
In particular, this implies   $c=1$, and $G'\cong A_n$. 

   Thus, we assume that $G'=A_n$. In this case  we (can) assume that $Z(G)=1$
and $G=A_n$ or $S_n$. (Indeed, let $z=g^{o(g)}$. If $z=1$ then $Z(G)=1$.
Otherwise, let $\phi(z)=\lam\cdot \Id$ and let $x=\phi(g)\zeta$, where $\zeta$ is  a primitive $o(g)$-root of $\lam\up$. Then $x^{o(g)}=1$  so $|x|=o(g)$. 
In addition, $\lan x,G'\ran$ is isomorphic to $A_n$ or $S_n$, and the \ei multiplicities
of $x$ and $g$ are the same as $x,g$ differ by a scalar multiple.) Then $G_1 =X_1\cong A_k$,   
and $G_3'=G_1\cong A_k$. We have shown above that 
$\tau|_{G_3'}$ is either subnatural or trivial. By Lemma \ref{kz2}, every \ir constituent of $\phi|_{A_n}$
is subnatural, whence the result.

 Suppose now that  
$o(g)\leq 4$.   Then we revise the choice of the  partition   $\Om=\Om_k\cup \Om_{n-k}$ as follows.  
If the cycle decomposition of $h$ has at least two cycles of order 4 then we take $k=8$ and we can assume that the restriction of $h$ to $\Om_k$ consists of two cycle of size 4. 
Otherwise we take $k=10$ and choose (a conjugate of) $g$ so  that the restriction of $h$ to $\Om_k$ be an even permutation of order $o(g)$. (Recall that $n>11$.)
Then $h$ is contained in a subgroup $Y$ of $G$ isomorphic to $A_k\times A_{n-k}$ or $A_k\times S_{n-k}$  if $h\in A_n$ and $h\notin A_n$, respectively. 
So  $h\in X_1\times X_2$, where $X_1\cong A_k$ 
and $X_2$ is isomorphic to $A_{n-k}$ or $S_{n-k}$. 

 Let $G_1$ 
be the preimage of $X_1$ in $G$  and $G_3=\lan g,G_1\ran$.   
Then $G_3'=c.A_k$ with $k\in\{8,10\}$, and $\tau(g)$ is almost cyclic for every \ir \ccc of $\phi|_{G_3}$. As in this case $k\leq 10$, we can use  Lemma \ref{a12} (instead of Lemma \ref{ct7} used above).  So either $\tau$ is almost cyclic or $\tau$ is as indicated in Table 2.1 or 2.2. However, none of the table entries for $k=8,10$ satisfies   $o(g)\leq 4$. This rules out the option $o(g)\leq 4$, and completes the proof.   \enp
   
\begin{lemma} \label{bp4} For  $n>11$ let G be a group such that $G/Z(G)\cong A_n$
 or $S_n$. Let $g\in G\setminus Z(G)$ be a p-element for $p>2$, and let $\phi: G\ra GL(m,F)$ be a faithful \irr such that    $\phi (g)$  is almost cyclic.   
Then $\phi(Z(G))=\Id$ and    $\phi|_{G'}$ are  subnatural.\el

\bp %Suppose first that either $p>2$ or $G/Z(G)\cong S_n$. 
As $\phi$ is faithful,  $Z(G)$ is a cyclic group, and then we can assume that  $G=G'$. (Indeed, $g\in G'Z(G)$ as $p>2$,  and if  $g=g_1z$ for $z\in Z(G),g_1 \in G'$ then $\phi(g)$ is almost cyclic \ii so is $\phi(g_1)$.) So we can assume that $G=c.A_n$ for $c\leq 2$, and then $G=G'$. 

Let $\Om$ be the natural permutation set for $A_n$ and let $h$ be the projection of $g$ into $G/Z(G)=A_n$. Let $k=o(g)$. By Lemma \ref{ct7}, we can assume that 
 $h$ stabilizes some subset of  $\Om_k\subset \Om$ of $k$
points,   $1\leq k<n$, and acts transitively on it.  So $h$ is contained in a subgroup $X\cong A_k\times A_{n-k}$  of $G/Z(G)$.
 Let $G_1 $ be the preimage of  $A_k$ in $G$,  so $G'\cong c.A_k$. 
Set $G_3=\lan g, G_1\ran$. Then $G_3'=c.A_k$ and $G_3/Z(G_3)\cong A_k$. 

 Suppose that $k>9$. The we claim that $c=1$ and the \ir constutuents of $\phi|_{G_3'}$ are trivial or subnatural. Suppose the contrary, and let 
$\tau$ be an \ir \ccc of $\phi|_{G_3}$ with $\dim\tau>1$. Then  $\tau(g)$ is almost cyclic. As $\Om_k$ is a maximal
$h$-orbit on $\Om$, we observe that $k=o(g)$ is the minimal integer $m$ such that $g^m\in Z(G_3)$.  By Lemma \ref{ct7}, $c=1$ and $\tau|_{G_3'}$ is 
subnatural. Whence the claim. As $c=1$, we have $G\cong A_n$.

\itf  every \ir \ccc of $\phi|{G_3'}$ is either trivial or subnatural. By Corollary \ref{kz2}, 
 $\phi$ is subnatural. So the result follows if $o(g)>9$.  

\def\cc{constituent}

%\med
 Suppose that $o(g)=|g|\leq 9$.  
We first fix some special cases.   Recall that $n>11$. For $n=12$  the result follows by Lemma \ref{a12}.  
If $n=13$ then $h$ fixes a point on $\Om$, so the result follows from the one with $n=12$. Let $n>13$. Suppose the contrary, that $\phi$ is not almost cyclic. If $n=14$ then $h$ fixes a point on $\Om$, unless $h$
is a double cycle of order 7. In the latter case, by Lemmas \ref{vw1} and \ref{va1}, we have  $\dim\phi\leq 12$. If $c=1$ then, arguing as in the proof of Lemma \ref{ct7}, using  \cite{J} we observe
that $\dim\phi\geq n(n-5)/2$, whence $n(n-5)/2\leq 4o(g)-4$; for $n=14$ this yields $63\leq 24 $, a contradiction. If $c=2$ then, by Lemma \ref{wo2}, $\dim\phi\geq 2^{5}=32$, which contradicts the above inequality $\dim\phi\leq 24$. 

The argument for $h$ a double cycle in $A_{14}$ works similarly for $h$ a double cycle of order 9 in $A_{18}$. So in this case   $\phi$ is almost cyclic. Note in addition, that, by Lemma \ref{33e}, $\phi(g)$ is not almost cyclic if  $g$ is a double cycle in 
$A_{2n}$ with $n$ odd.    

 In general, let $n>14$.    
Then we revise the above partition   $\Om=\Om_k\cup \Om_{n-k}$, by choosing 
  $k=9,10,14$ for $o(g)=3,5,7$, respectively, and $k=12,18$ if $o(g)=9$ and a cycle  of order 9 occurs once or twice, respectively, in the cycle decomposition of $h$. We mimic the above reasoning to deduce that $c=1$ and $\phi$ is subnatural. 

Let  $\Om_k,G_1,G_3=\lan g,G_1\ran$  be as above for $k$ just specified.  
By the above, we are left with the cases where the cycles of $h$ on $\Om_1$ are of maximal possible sizes. So $(o(g),k)\in\{(3,9), (5,10), (7,14),(9,12)$, $(9,18)\}$.
Recall that $G_3/Z(G_3)\cong A_k$ and $G_3'\cong c.A_k$ for $c\leq 2$. Let $\tau$ be an \ir \ccc of $\phi|_{G_3}$ with $\dim\tau>1$. Then $\tau(g)$
is almost cyclic. We claim that $\tau|_{G_3'}$ is subnatural. If $k\leq 12$ then the claim follows from Lemma \ref{a12} and Table 1,
if $k=14,18$ then this is proved in the previous paragraph. So $c=1$, $G\cong A_n,$ $G_3'\cong A_k$ and $\tau|_{G_3'}$ is \irt 

Therefore, every \ir \ccc of $\phi|_{G_3'}$ is either trivial or subnatural. So the result follows from Corollary \ref{kz2}.\enp

\begin{corol}\label{vv6} Let $G=2.A_n$, $n>4$, 
and let $g\in G$ be a p-element, $p>2$. Let $\phi:G\ra GL_m(F)$ be a non-trivial   \ir   \rep of G.
% and $g\in G/Z(G)$ with $| prime power. 
Suppose that $\ell\neq p$ and  $\phi(g)$ is almost cyclic. 
If $\deg\phi(g)< o(g)-1$ then $o(g)\leq 9$ and one of the \f holds:
\med

$(1)$  $|g|=9$,  $m=8$, $\deg\phi(g)=7$
and either $n=9$,   $\ell\neq 2$ or  $n=10$,  $\ell=5;$  

\med
$(2)$ $|g|=5$,  $\deg\phi(g)= 3$, $m=3$ and either  $n=5$,   $\ell\neq 2$ or  $n=6$, $\ell=3;$

\med
$(3)$   $|g|=7$,  $m=4$, $\deg\phi(g)= 4$ and either $n=7$,  or $n=8$,  $\ell=2;$ 

\med
$(4)$ $|g|=5$,  $m=2$ and $\deg\phi(g)= 2$ and either  $n=6$, $\ell=3$ or $n=5.$  

%\med $(5)$ $|g|=7$,  $n=7$, $\ell=5$,   $m=3$ and $\deg\phi(g)=3$. 
\end{corol}  

\bp If $\phi$ is subnatural then the result follows from Lemma \ref{33e}. Otherwise, $n\leq 11$ by  Lemmas \ref{sc0} and \ref{bp4}. For $n\leq 11$
the result follows by Lemma \ref{a12}  and Table 1.\enp

Remark. If $G=c.A_n$, $c>2$ and  $\deg\phi(g)< o(g)-1$ then  we additionally have $G=3A_7$, $|g|=7$ and $(m, \deg\phi(g),\ell)=(3,3,5)$ or $(4, 4,\neq 2,)$,
%$\ell=5$, $m=3$ and  $\deg\phi(g)\geq o(g)-4$, 
or $G=3A_6$, $|g|=5$ and $(m, \deg\phi(g),\ell)=(3,3,5)$.

\begin{corol}\label{v77} Let $G=2.A_n$, or $2.S_n$, $n>4$, 
and let $g\in G$ be a $2$-element. Let $\phi:G\ra GL_m(F)$ be a non-trivial   \ir   \rep of G.   Suppose that $\ell\neq 2$ and  $\phi(g)$ is almost cyclic. 
If $\deg\phi(g)< o(g)-1$ then $o(g)=4$, $m=2$ and either

$(1)$ $n=6$,  $G=2.A_6$,   $\ell=3$ and $\deg\phi(g)=2$ or

$(2)$ $n=5$,  $G=2.S_5$,   $\ell=5$  and $\deg\phi(g)= 2.$\end{corol}

\bp Repeat the reasoning in the proof of Corollary  \ref{vv6} with use  Lemma \ref{sc0}   instead Lemma \ref{bp4}. \enp

If, more generally, we have  $G' \cong c.A_n$, $c>2$, then we detect two more cases:     $G=3.A_6.2_2$,  $\ell\neq 2$,  $o(g)=8$, %or $G=2.A_6.2_3$, $\ell=3$, 
$m=\deg\phi(g)=3$ and      $G=3.A_6.2_3$, $o(g)=8$,  $\ell\neq 2,3$, $m=\deg\phi(g)=6.$

\def\irc{irreducible constituent }

\medskip
\bp[Proof of Theorem {\rm \ref{th1}}] If $n\leq 11$, the result follows by
inspection of  Tables 1,1 and 2.1. If $n>11$ then Lemma \ref{sc0} for $p=2$ and Lemma \ref{bp4} for $p>2$ imply the result.  
\enp

%COMPARE WITH CRAVEN p.15, Pr 4.1. In fact, if $\phi\mod p$ is \ir and g is AC then $ g\pmod{p}$ is AC, so there is some points to check.

 Department of Mathematics, University of Brasilia, Brasilia-DF 70910-900, Brazil 

 e-mail: alexandre.zalesski@gmail.com

\newpage
\begin{center}Table 1.1. The non-trivial {non-natural } irreducible $\ell$-modular \reps $\phi$ \\ of $G=c.A_n$, $n\leq 10$, and $p$-elements $g\in G$, $p>2$, with $\phi(g)$  almost cyclic
 % in non-subnatural \ir \reps $\phi$ of $G=c.A_n$, $n\leq 12$

\bigskip

\small{
 \begin{tabular}{ |c|c|c|c|c|c|c|c|}
\hline
$n$&$\ell$&$\dim\phi$&$c$&$o(g)$&$\deg\phi(g)$&$e$&$m$\\
\hline 
$5$&$2$&$2$&$1$&$3,5$&$2$&$-$&$1$\cr
\hline
$5$&$\neq 2,p$&$2$&$2$&$3,5$&$2$&$-$&$1$\cr
\hline 
$5$&$\neq 2,p$&$3$&$1$&$3,5$&$3$&$-$&$1$\cr
\hline 
%$5$&$2$&$4 sbn?$&$1$&$3$&$3$&$1$&$2$\cr
%\hline 
%$5$&$2$&$4 snb?$&$1$&$5$&$4$&$-$&$1$\cr
%\hline 
$5$&$\neq 2,3$&$4$&$2$&$3$&$4$&$1$&$2$\cr
\hline 
$5$&$\neq 2,5$&$4$&$2$&$5$&$4$&$-$&$1$\cr
\hline 
$5$&$\neq 2,3,5$&$5$&$1$&$5$&$5$&$-$&$1$\cr
\hline 
$5$&$\neq 2,5$&$6$&$2$&$5$&$4$&$1$&$2$\cr
\hline 
\hline
$6$&$3$&$2$&$2$&$5$&$2$&$-$&$1$\\
\hline
%$6$&$\neq 2,3$&$3$&$3$&$2$&$2$&$-1$&$2$\\
%\hline
%$6$&$3$&$3$&$1$&$2$&$2$&$-1$&$2$\\

%\hline 
$6$&$ \neq 3,p$&$3$&$3$&$3,5$&$3$&$-$&$1$\\
\hline 

$6$&$3$&$3$&$1$&$5$&$3$&$-$&$1$\\
%\hline 
%$6$&$2$&$4 sbn$&$1$&$3B$&$3$&$1$&$2$\\
%\hline
%$6$&$ 2$&$4 sbn$&$1$&$5$&$4$&$-$&$1$\\
\hline
$6$&$ 2$&$4$&$1$&$3A$&$3$&$1$&$2$\\
\hline
$6$&$ 2$&$4$&$1$&$5$&$4$&$-$&$1$\\
\hline
$6$&$\neq 2,3$&$4$&$2$&$3B$&$4$&$1$&$2$\\
\hline
$6$&$\neq 2,5$&$4$&$2$&$5$&$4$&$-$&$1$\\
\hline
$6$&$\neq 2,3,5$&$5$&$1$&$5$&$5$&$-$&$1$\\
\hline
$6$&$\neq 2,3$&$5$&$1$&$3B$&$3$&$1$&$3$ \\
\hline
$6$&$3$&$6$&$2$&$5$&$5$&$1$&$2$\\
\hline
$6$&$>5$&$6$&$3,6$&$5$&$5$&$1$&$2$\\
\hline
\hline 
$7$&$5$&$3$&$3$&$3,7$&$3$&$-$&1\\
\hline  
$7$&$2$&$4$&$1$&$3B$&$3$&$1$&$2$\\
\hline  
$7$&$2$&$4$&$1$&$5,7$&$4$&$-$&$1$\\
\hline 
$7$&$\neq 2,3$&$4$&$2$&$3B$&$3$&$1$&$2$\\
\hline
$7$&$\neq 2,p$&$4$&$2$&$5,7$&$4$&$-$&$1$\\
\hline
%\hline
$7$&$\neq 2,3,p$&$6$&$3,6$&$5$&$5$&$1$&$2$\\
\hline
$7$&$\neq 2,3,7$&$6$&$3,6$&$7$&$6$&$-$&$1$\\
\hline
$7$&$ 2$&$6$&$3$&$5$&$5$&$1$&$2$\\
\hline
$7$&$ 3$&$6$&$2$&$7$&$6$&$-$&$1$\\
\hline
\hline
$8$&$2$&$4$&$1$&$3B$&$3$&$1$&$2$\\
\hline
$8$&$2$&$4$&$1$&$5,7$&$4$&$-$&$1$\\
\hline
%$8$&$2$&$6$&$1$&$5$&$5$&$1$&$2$\\\hline
%$8$&$2$&$6$&$1$&$7$&$6$&$-$&$1$\\\hline 
%$8$&$\neq  2$&$8$&$2$&$8$&$8$&$-$&$1$\\
%\hline 
$8$&$\neq  2,7$&$8$&$2$&$7$&$7$&$1$&$2$\\
\hline
\hline 
$9$&$2$&$8$&$1$&$7$&$7$&$1$&$2$\\
\hline 
$9$&$\neq  2,7$&$8$&$2$&$7$&$7$&$1$&$2$\\
\hline 
$9$&$2$&$8$&$1$&$9$&$8$&$-$&$1$\\
\hline 
$9$&$2$&$8$&$1$&$9$&$7$&$1$&$2$\\
\hline
$9$&$\neq2,3$&$8$&$2$&$9A$&$7$&$1$&$2$\\
\hline
$9$&$\neq2,3$&$8$&$2$&$9B$&$8$&$-$&$1$\\
\hline
\hline 
%$9nat?$&$ 3$&$7$&$1$&$7$&$7$&$-$&$1$\\\hline 
$10$&$5$&$8$&$2$&$7$&$7$&$1$&$2$\\
\hline
$10$&$5$&$8$&$2$&$9A$&$7$&$1$&$2$\\
\hline
$10$&$5$&$8$&$2$&$9B$&$8$&$-1$&$2$\\
\hline
 
%$S_5$&$\neq 2,5$&$4$&$1?$&$5$&$4$&$-$&$1$\\\hline

\end{tabular}
}
\end{center}

In the tables  $e$ stands for the \ei of $\phi(g)$ of \mult greater than 1, if it exists, and $m$ is  the maximal \ei \mult of $\phi(g)$. See comments after the proof of Lemma \ref{} 0000000 on the last two lines in Table 2.1.  In column $o(g)$ we use notation of \cite{AB} to specify the conjugacy class of a certain order elements if this is necessary.  
%he entries in the $e$-column and the last  row of Table 2.1 we write $\pm 1$ to say that in two \ir \reps of $2.A_n$ of degree 8the eigenvalues of \mult 2 of elements $g\in G$ of order 8   differ by sign.  

\newpage

\begin{center}Table $1.2$. Almost cyclic elements $g\in G'$, $|g|$ a $p$-power, $p>2$,
  in non-trivial \\ {\it non-subnatural} \ir \reps of $G$ with $G/Z(G)=S_n$, $4<n< 13$

\bigskip

\small{
 \begin{tabular}{ |c|c|c|c|c|c|c|c|}
\hline
$G/Z(G)$&$\ell$&$\dim\phi$&$|Z(G')|$&$o(g)$&$\deg\phi(g)$&$e$&$m$\\
\hline 
$S_5$&$\neq 2,5$&$6$&$1$&$5$&$4$&$1$&$2$\\
\hline
$S_5$&$\neq 2,5$&$4$&$1$&$5$&$4$&$-$&$1$\\
\hline
$S_6,A_6.2_3$&$3$&$4$&$2$&$5$&$4$&$-$&$1$\\
\hline
$S_6,A_6.2_3$&$\neq 3$&$6$&$3$&$5$&$4$&$1$&$2$\\
\hline
$S_6,A_6.2_3$&$ 3$&$6$&$1$&$5$&$4$&$1$&$2$\\
\hline
$S_7$&$\neq 7$&$8$&$2$&$7$&$6$&$1$&$2$\\
\hline
$S_8$&$2$&$8$&$1$&$7$&$6$&$1$&$2$\\
\hline
 
\end{tabular}
}
\end{center}

%\newpage

%In Table 2 $c.A_n$ or $c.S_n$ means a central extension of $A_n,S_n$, respectively, such that $c.A_n$ is perfect, and the derived subgroup of  $c.S_n$ is $c.A_n$. 

%In addition, $2.A_6.2_1$ is a group whose derived subgroup  is $c.A_6$ and $2.A_6.2/Z(2.A_6.2)$ is not isomorphic to $S_6$.
%$2.A_6.2_2\cong PGL_2(9)$, and $2.A_6.2_3\cong PGL_2(9).2\not\cong PGL_2(9)$. Observe that there are two conjugacy classes of elements of order 9 in $A_9$ and $A_{10}$,
%and two  \ir characters of degree 8 which are  not subnatural.

\bigskip
\centerline{Table 2.1. Almost cyclic elements of $2$-power order in faithful \ir }

\centerline{\reps of $G$, where $G/Z(G)=A_n$, $n<12$} %and $Z(G)=Z(G')$ } % and their  central extensions ($S_6$ to be refined)}

\bigskip
\begin{center}

 \begin{tabular}{ |c|c|c|c|c|c|c|}

\hline
$G'$&$\ell$&$\dim\phi$&$o(g)$&$\deg\phi(g)$&$e$&$m$\\

\hline 
$2.A_5$&$\neq 2$ &$2$&$2$&$2$&$ -$&$1$\\
\hline 
$A_5$&$\neq 2$ &$3$&$2$&$2$&$ -1$&$2$\\
\hline 
$2.A_6$&$3$ &$2$&$2,4$&$2$&$ -$&$1$\\
\hline 
$A_6$&$3$ &$3$&$2$&$2$&$ -1$&$2$\\
\hline 
$A_6$&$3$ &$3$&$4$&$3$&$ -$&$1$\\
\hline 
$3.A_6$&$\neq 2,3$ &$3$&$2$&2&$ -1$&$2$\\
\hline 
$3.A_6$&$\neq 2,3$ &$3$&$4$&$3$&$ -$&$1$\\
\hline 
$3.A_7$&$5$ &$3$& $2$&$2$&$ -1$&$2$\\
\hline 
$2.A_7$&$\neq 2$ &$4$&$4$&$4$&$ -$&$1$\\
\hline 
$2.A_7$&$\neq 2$ &$4$& $4$&$4$&$ -$&$1$\\
\hline
$2.A_{10}$&$5$ &$8$& $8$&$7$&$ \pm 1$&$2$\\
\hline
$2.A_{10}$&$5$ &$8$& $8$&$7$&$\mp  1$&$2$\\
\hline
\end{tabular}
\end{center}

%\newpage 
\bigskip
\centerline{Table 2.2. Almost cyclic elements of $2$-power order in {\it non-subnatural}     faithful }

\centerline{\ir \reps of $G=c.A_n.2$,  $n<12$, and $G=\lan g,G'\ran$} %and $Z(G)=Z(G')$ } % and their  central extensions ($S_6$ to be refined)}
%where $G/Z(G)\subset {\rm Aut}\,(A_n)$,
  
\bigskip
\begin{center}
 \begin{tabular}{ |c|c|c|c|c|c|c|}
\hline
$\lan g,G'\ran$&$\ell$&$\dim\phi$&$o(g)$&$\deg\phi(g)$&$e$&$m$\\

\hline 
$2.S_5$&$\neq 2$ &$2$ &$2,4$&$2$&$ -$&$1$\\
\hline 
%$S_5$&$\neq 2,3$ &$3$&$1$&$2$&$2$&$ -1$&$2$\\
%\hline 
%$S_5$&$\neq 2,3,5$ &$3$&$1$&$4$&$3$&$ -$&$1$\\
%\hline 
%$S_5$&$ \neq 2,3$ &$4? sbn$&$1$&$4$&$4$&$ -$&$1$\\
%\hline 
%$2.S_5$&$ 3$ &$4$&$2$&$4$&$4$&$ -$&$1$\\
%\hline 
%$2.S_5$&$ 5$ &$4$&$2$&$4$&$4$&$ -$&$1$\\
%\hline 
$2.S_5$&$\neq 2$ &$4$& $4$&$4$&$ -$&$1$\\
\hline 
$S_5$&$\neq 2,3 $ &$5$&$4$&$3$&$\pm 1$&$2$\\
%\hline 
%$2.S_5$&$ \neq 2,3$ &$6$&$2$&$5$&$5$&$ 1$&$2$\\

\hline
$S_6$&$3$ &$3$&$8$&$3$&$ -$&$1$\\
\hline 
$3.S_6$&$\neq 2,3$ &$3$&$4,8$&$3$&$ -$&$1$\\
\hline 
$2.A_6$&$\neq 2,3$ &$4$& $4$&$4$&$ -$&$1$\\
\hline 
$2.A_6.2_1$&$5$ &$4$& $4$&$4$&$ -$&$1$\\
\hline
$A_{6}.2_3$&$3$ &$6$& $8$&$6$&$-$&$1$\\
\hline 
$3.A_6.2_3$&$\neq 2,3$ &$6$ &$8$&$6$&$ -$&$1$\\
\hline 
$A_6.2_3 $&$5$ &$8$& $8$&$8$&$ \pm1$&$2$\\
\hline 
$S_6$&$\neq 2,3$ &$8$& $8$&$8$&$ -$&$1$ \\
\hline 
$2.S_6$&$\neq 2,3$ &$8$ &$8$&$8$&$ -$&$1$\\
\hline 
$2.S_8$&$\neq 2$ &$8$& $8$&$7$&$\pm \sqrt{-1}$&$2$\\
\hline 
%$2.S_8$&$\neq 2$ &$8$& $8$&$7$&$-\sqrt{-1}$&$2$\\
%\hline 
%$S_9$&$\neq 2$ &$8$&$1$&$8$&$7$&$1$&$2$\\
%\hline 
%$S_9$&$\neq 2$ &$8$&$1$&$8$&$8$&$-1$&$2$\\
% \hline 
$2.S_9$&$3$ &$8$& $8$&$7$&$\pm \sqrt{-1}$&$2$\\
\hline 
\end{tabular}
\end{center}

\end{document}